\documentstyle[11pt]{article}
\begin{document}
\title{ {Ideal Class Groups and Subgroups of Real Quadratic\\  Function
Fields}
         \thanks{MR (1991) Subject Classification: 11R58; 11R29; 14H05.}
         \thanks{Project Supported by the NNSFC (No.19771052)}
         \thanks{Email: kpwang@mail.cic.tsinghua.edu.cn}}
\author{ {WANG Kunpeng  $\; \; $ and $\; \; $ ZHANG Xianke }\\
         \smallskip \small{Department of Mathematical Sciences, Tsinghua
University, Beijing 100084, P.R. China }
       }
\date{}
\maketitle
\parindent 15pt
\baselineskip 20pt
\parskip 10pt

In [4], algebraic  quadratic number fields $F={\bf Q}(\sqrt{d})$ were studied,
 a necessary and sufficient
condition was obtained
 for the ideal class group $H(F)$ of any real field $F={\bf Q}(\sqrt{d})$  to have
  a cyclic subgroup of order $n$;
  and eight series of such fields $F$ were given explicitly
  by utilizing the theory of continued fractions. All these fields $F$
   are ERD-type or GERD-type,
that is, the discriminates $d$ are in the form $\; d=N^2+n,\; $
 where $\; n|4N\; $ or $\; n|2^m N \; \;
(m\geq 2).\; $ In the present note, we study algebraic function fields $K$,
 give necessary and sufficient
 condition for the ideal class group $H(K)$ of any real quadratic function
 field $K$ to have a cyclic subgroup of order $n$, and obtain
  eight series  of such fields $K$, with four of them not
  ERD-type or GERD-type.

Now suppose that $k={\bf F}_q(T)$ is the rational function field
 of indeterminate (variable) $T$ over ${\bf F}_q$, the finite field with $q$
 elements. Let $R={\bf F}_q[T]$ be
   polynomials over ${\bf F}_{q}$ , which is said to be the
  ring (domain) of integers (integral functions) of $k$.
Any finite algebraic extension $K$ of $k$ is said to
 be an algebraic function field.
 The study of algebraic function fields is equivalent
 to the arithmetic of algebraic smooth curves over ${\bf F}_q$.
On the other hand, the study of function fields is parallel to the study of
 algebraic number fields (via valuation theory).

  If  further  $r|f$, then $K$ and $D$ are said to be of ERD-type.
  The integral closure of $R$ in any algebraic function field
   $K$ is said to be the ring (domain) of
  integers of $K$, and is denoted by $ {{\cal O}}_K$.
  $ {{\cal O}}_K$
  is a Dedekind domain, and $K$ is the quotient field of
  $ {{\cal O}}_K$.
  The fractional ideals of $  {{\cal O}}_K$ form a multiplication group
  $  {{\cal I}}_K $. Let $ {{\cal P}}_K$ denotes
  the principal ideals in
  $  {{\cal I}}_K $. Then the quotient group
   $H(K) = H(D) = {{\cal I}}_K /  {{\cal P}}_K $
    is said to be the ideal
   class group of $K$.  And  $h(K) = h(D) = \#H(K)$
   (the order of $H(K)$) is said to be the
   ideal class number of $K$.
The theory of continued fractions plays a vital role in the research of real quadratic
number fields.
 We here will apply the
continued fraction theory in [3]  to study real quadratic function fields,
 obtaining our main results on their ideal class groups and subgroups.

{\bf Theorem 1.} The ideal class group  $H(D)$ of a real quadratic function field
$K=k(\sqrt{D})$ contains a subgroup of order $n\ (\geq 2)$ if and only if the equation
$$X^2-DY^2=cZ^n \; \; \quad  (Z  \in R-{\bf F}_q,\;\; c\in {\bf F}_q^\times )$$
 has a solution
$(X,Y)$ where $X,Y\in R$ are relatively prime,
 and equations $$X^2-DY^2=c'Z^j \; \; \;  \quad
(c'\in {\bf F}_q^\times, \; \; 1\leq j|n,\;  j\not=n)$$ have no solution with $(X,Y)=1$.

If $(X,Y)\in R\times R$  is a solution of the equation
 $X^2-DY^2=C\; (C\in R)$, then we say $\alpha =X+Y\sqrt{D}$ is a solution.
 If  $(X,Y)=1$ then $\alpha $ is said to be a primary solution.
  For any unit $\varepsilon $  of the Dekekind domain ${\cal O}_K=R+R[\sqrt{D}]$,
  $\varepsilon \alpha $ is said to be an associate of  $\alpha $.
   The conjugate of   $\alpha =X+Y\sqrt{D}$ is $  \bar{\alpha } =
   \sigma ( \alpha ) =X-Y\sqrt{D}$, where $\mbox{Gal}(K/k)=\langle \sigma
\rangle$.

{\bf Lemma.} The primary solutions of the equation  $X^2-DY^2=C\; $  $\; (C\in R, \;
 \mbox{deg}C<\frac{1}{2}\mbox{deg}D)$  are just
$$p_{i-1}+q_{i-1}\sqrt{D}\; \qquad \; ( (-1)^i Q_i=C,\; 0<i\in {\bf Z}) $$
 together with their associates and
conjugates ,  where $p_{i-1}/q_{i-1}=[a_0,a_1,\cdots,a_{i-1}]$ is the $(i-1)-$th
convergent of the simple continued fraction
 $\sqrt{D}=[a_0,a_1,\cdots]$; $\; Q_i$ is the denominator of
the  $i$-th complete quotient $\alpha_i$ $=[a_n,a_{n+1},\cdots]$ $=(\sqrt{D}+P_i)/Q_i$.

In the following Theorem 2 and 3, we give explicitly eight series
 of real quadratic
function fields whose class group contains a cyclic subgroup of order $n$. By the above
Lemma and Theorem 1, we could use the simple continued fraction of $\sqrt{D} $
 to prove the class group of the field $K=k(\sqrt{D})$ to contain
  a  cyclic subgroup of order $n$. If $K$ is  a field of
 ERD-type, then it is easy to find the expansion of the simple
 continued fraction of $\sqrt{D}$.
 The fields in Theorem 2 are in this case, i.e., they are ERD-type.
   Fortunately, for some other $\sqrt{D}$ not being ERD-type,
   we could also obtain the expansions of
   simple continued fractions,
   which gives Theorem 3.

{\bf Theorem 2.} For the following series of monic square-free polynomials
 $D\in R$, the ideal class groups $H(D)$ of the real quadratic function field
  $K=k(\sqrt{D})$  all contain a subgroup of order $n$ (Here we assume
   $ Z\in R$ is any non-constant polynomial). \\
\mbox{}\hspace{1cm} (1) $D=Z^{2n}+1;$ \\ \mbox{}\hspace{1cm} (2) $D=(Z^n+F-1)^2+4F,$
\mbox{}\hspace{0.65cm} where $F|(Z^n-1); $\\ \mbox{}\hspace{1cm} (3) $D=(Z^n-F+1)^2+4F,$
\mbox{}\hspace{0.65cm} where $F|(Z^n+1); $
\\ \mbox{}\hspace{1cm} (4) $D=(4Z^n+F-1)^2+4F,$ \mbox{}\hspace{0.5cm}
where $F|(4Z^n-1).$\\ In general, for any $$D=(Z^n+aF-b)^2+4abF, $$
 where  $a,b\in {\bf F}_q^\times,\;\;
F|(Z^n-b),\; $ the ideal class of  $K=k(\sqrt{D})$
 contains a subgroup of
order $n$.

{\bf Theorem 3.} Suppose that $a\in R-{\bf F}_q$ is a non-constant polynomial,
 the polynomial $ A=2a+1\in R$ is monic. Let   $A=Z^n$ for any non-constant
  $Z\in R-{\bf F}_q$. Then for the following square-free polynomials $D$,
the ideal class groups $H(D)$ of $K=k(\sqrt{D})$ all contain a cyclic subgroup
 of order $n$. \\
\mbox{}\hspace{1cm}(1) $D=(A^d+a)^2+A;$ \\ \mbox{}\hspace{1cm}(2) $D=(A^d-a)^2+A;$\\
\mbox{}\hspace{1cm}(3) $D=(A^d+a+1)^2-A;$ \\ \mbox{}\hspace{1cm}(4) $D=(A^d-a-1)^2-A.$

{\bf Corollary.} Suppose that the polynomials $D$ are as in theorem 2 and 3,  then
 the ideal class numbers $h(D)$ of  $K=k(\sqrt{D})$ all have a factor $n$.
 in particular, we have $h(D)\geq n$.
                 \vskip 1cm

{\bf References:}
 \baselineskip 0pt
\parskip 0pt
\begin{description}

\item[[1]] E. Artin, Quadratische K\"orper im Gebiete der h\"oheren
Kongruenzen I, II,
          Math. Z., {\bf 19} (1924), 153-206, 207-246.
\item[[2]] S. Louboutin, Continued fraction and real quadratic fields, J.
Number Theory,
        {\bf 30} (1998), 167-176.
\item[[3]] WANG Kunpeng, ZHANG Xianke , Continued fractions of
functions connected with real quadratic function
         fields (to appear)
\item[[4]] ZHANG Xianke, L. C.  Washington, Ideal class groups and there
subgroups of real quadratic fields, Sci. in China,  (A) {\bf 27}  (1997), 522-528.

\end{description}

\end{document}